\documentclass{article}
\usepackage{amsmath}
\usepackage{amsfonts}
\usepackage{amsthm}
\newcommand{\RR}{\mathbb{R}}
\newtheorem{theorem}{Theorem}
\newtheorem*{thm1}{Theorem 1}
\newtheorem{lemma}{Lemma}
\theoremstyle{definition}
\newtheorem{defn}{Definition}
\theoremstyle{remark}
\newtheorem*{acn}{Acnowledgements}
\DeclareMathOperator{\orn}{or}
\DeclareMathOperator{\inv}{inv}
\begin{document}
\title{Tensor invariants of $SL(n)$, wave graphs and L-tris 
\thanks{This research was partially supported by 
ONR grant N00014-97-1-0505.}}

\author{Aleksandrs Mihailovs\\
Department of Mathematics\\
University of Pennsylvania\\
Philadelphia, PA 19104-6395\\
mihailov@math.upenn.edu\\
http://www.math.upenn.edu/$\sim$mihailov/
}
\date{\today}
\maketitle

\begin{abstract}
The space of invariants of a tensor product of representations of $SL(n)$ 
is provided with the basis parametrized by wave graphs introduced here 
especially for this purpose. The proof utilizes a game similar to Tetris,  
named here L-tris.
\end{abstract}    
\setlength{\baselineskip}{1.5\baselineskip}

\section{Introduction}
In my earlier work \cite{M1,M2} among other results, 
I gave the new proof of the famous 
theorem of Rumer, Teller and Weyl \cite{W} parametrizing a basis of 
the subspace of $SL(2)$-invariants of 
$V^{\otimes m}=V \otimes \dots \otimes V$ ($m$ times), 
where $V$ is the two-dimensional linear space 
with the standard action of $SL(2)$, by 
$1$-regular graphs with the vertices $1,2, \dots, m$, 
edges of which can be drawn in the upper half-plane without intersections. 

In this work I construct a basis of the subspace of $SL(n)$-invariants of 
$V^{\otimes m}$, where $V$ is the $n$-dimensional linear space with the 
standard action of $SL(n)$, by the $n$-wave graphs with $m$ vertices. 
An $n$-wave graph is a graph with the vertices $1,2, \dots, m$, each 
connected component of which is a path of length $n-1$, 
edges of which  
can be drawn in the book with $n-1$ pages, i.\ e.\ 
$n-1$ copies of the upper half-plane, glued along $\RR$, such that the 
first edge of each connected component, $\{i_1i_2\}$, is drawn on the first page, the 
second edge, $\{i_2i_3\}$, on the second page, and so on, such 
that $i_1<i_2<\dots<i_n$, without intersections. 

$2$-wave graphs are exactly the same graphs as were used for the 
description of $SL(2)$-invariants. 
Here are two of a total number of five of $3$-wave graphs with $6$ vertices:

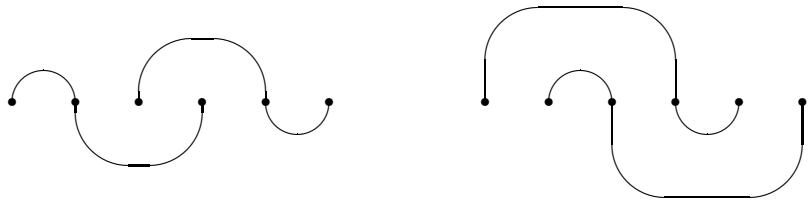
\begin{picture}(410, 84) 
  
\put(33,48){\circle*{3}}
\put(57,48){\circle*{3}}
\put(81,48){\circle*{3}}
\put(105,48){\circle*{3}}
\put(129,48){\circle*{3}}
\put(153,48){\circle*{3}}
\put(45,48){\oval(24,24)[t]}
\put(81,48){\oval(48,48)[b]}
\put(105,48){\oval(48,48)[t]}
\put(141,48){\oval(24,24)[b]}

\put(212,48){\circle*{3}}
\put(236,48){\circle*{3}}
\put(260,48){\circle*{3}}
\put(284,48){\circle*{3}}
\put(308,48){\circle*{3}}
\put(332,48){\circle*{3}}
\put(248,48){\oval(72,72)[t]}
\put(248,48){\oval(24,24)[t]}
\put(296,48){\oval(72,72)[b]}
\put(296,48){\oval(24,24)[b]}

\end{picture}

The corresponding invariants are 
\begin{equation}\label{1}
\left( (x\wedge y\wedge z)\otimes(x\wedge y\wedge z)\right)^{\sigma_i},\quad i=1,2,
\end{equation}
where
\begin{equation}\label{2}
\sigma_1=(34)\in S_6,\quad \sigma_2=(24)(35)\in S_6.
\end{equation}
Note that the choice \eqref{2} of the permutations $\sigma_i$ in \eqref{1}, is not 
unique for the given wave graphs. Multiplying $\sigma_i$ by $(14)(25)(36)$, we 
obtain the same result.

For a wave graph $G$, denote $t_G$ the analogous tensor products of the basic 
wedge $n$-forms, see \eqref{10}. 

\begin{thm1}
Tensors $t_G$ parametrized by all $n$-wave graphs with $m$ vertices, form a basis 
in the space of $SL(n)$-invariants in $V^{\otimes m}$.
\end{thm1}

The proof uses a game similar to Tetris, named here L-tris.

\section{The main theorem}
In this section we give all the necessary definitions and prove the main theorem.

Let $f$ be a field of characteristic $0$, and $SL(n)$ ---the group of 
$n \times n$ $f$-matrices with determinant 1, acting on $n$-dimensional 
linear $f$-space $V$ with basis $(x_1, \dots, x_n)$ by the standard way.

Recall some fundamental facts about the 
representations of $SL(n)$. The word {\em representation} will mean below 
a polynomial finite dimensional linear representation over $f$.
Every representation 
of $SL(n)$ is equivalent to a sum of irreducible representations. 
All classes of equivalence of the irreducible 
representations are parametrized by partitions of length $<n$. Denote 
$L_n$ the set of partitions of length $<n$, and
denote $\rho_\lambda$ the irreducible representation corresponding to a 
partition $\lambda\in L_n$. 
Then $\rho_0$ is a trivial representation of dimension $1$ 
and $\rho_1$ is the standard representation in $V$ mentioned above.
To describe the decomposition of tensor products 
of irreducible representations, we'll use Young diagrams.

For a moment, let us draw the Young diagrams rotated by $90^\circ$ counterclockwise.
Then we can interpret a Young diagram of a partition of length $<n$ as a 
Tetris position on a Tetris game field of width $n$, with non-increasing height 
of columns (looking from the left hand side to the right). 

\begin{defn}
For a partition $\mu$ of length $<n$, 
denote $T_n(\mu)$ the set of partitions, Young diagram 
of which can be obtained from the Young diagram of $\mu$ by dropping to it a 
$1\times 1$ block and applying the Tetris rules, i.\ e.\ deleting of rows of 
length $n$ if possible.
\end{defn}

Then
\begin{equation}\label{5}
\rho_\mu\otimes\rho_1\simeq\sum_{\lambda\in T_n(\mu)}\rho_\lambda
\end{equation}

\begin{lemma}\label{lem1}
\begin{equation}\label{6}
\rho_1^{\otimes m}\simeq\sum_{\substack{\lambda\in L_n, |\lambda|\leq m\\
|\lambda|\equiv m \bmod n}}f^{\tau(\lambda,m,n)} \rho_\lambda,
\end{equation}
where $|\lambda|$ denotes the weight (i.\ e. the sum of all parts) of 
a partition $\lambda$; $\tau(\lambda,m,n)=\lambda\bigcup n^{(m-|\lambda|)/n}$
is the partition of weight $m$,  
Young diagram of which can be obtained from the Young diagram of $\lambda$ 
by adding the necessary number of Tetris rows of length $n$,
and $f^\tau$ is the number of standard tableaux of shape $\tau$.
\end{lemma}
\begin{proof} By induction on $m$, by iteration of \eqref{5}, 
the Young diagrams of the partitions $\lambda$ in the right hand side of \eqref{6} 
can be obtained by dropping $m$ $1\times1$ blocks on the Tetris game field 
and applying the Tetris rules. Sequential numbering of dropping blocks from $1$ 
to $m$ gives us the standard tableaux of shape $\tau$. Conversely, each standard 
tableaux defines the places where to put each dropping block. Thus, we have a 
bijection between the set of standard tableaux of shape $\tau(\lambda,m,n)$ and 
the set of irreducible components of type $\rho_\lambda$ in the total 
decomposition of the left hand side of \eqref{6} according to \eqref{5}. 
\end{proof} 
 
\begin{lemma}\label{lem2}
The subspace of $SL(n)$-invariants 
of $V^{\otimes n}$ is one-dimensional, and we can choose 
\begin{equation}\label{3}
\omega=x_1 \wedge x_2 \wedge \dots \wedge x_n
\end{equation} 
as a basis element in that space.
\end{lemma}
\begin{proof}
For $a\in SL(n)$, 
\begin{equation}\label{4}
a \omega=\det a\cdot\omega=\omega
\end{equation}
It means that $\omega$ is invariant.
By Lemma \ref{lem1}, the dimension of the space of invariants is equal to 
$f^{\tau(0,n,n)}=f^n=1$. 
\end{proof}

\begin{lemma}\label{lem3}
The dimension of the subspace of $SL(n)$-invariants of $V^{\otimes m}$ is nonzero 
iff $m=kn$ for some integer $k$, in which case it equals  
\begin{equation}\label{7}
\frac{m!\thinspace (1!\thinspace 2! \dots (n-1)!)}{k!\thinspace (k+1)! \dots (k+n-1)!} .
\end{equation}
\end{lemma}
\begin{proof}
It follows from Lemma \ref{lem1} that the dimension equals $f^{n^k}$. Using 
the hook length formula \cite{FRT}, one obtains \eqref{7}.
\end{proof}  

\begin{defn}
An $n$-wave graph is a graph with the vertices $1,2, \dots, m$, each 
connected component of which is a path of length $n-1$, edges of which  
can be drawn in the book with $n-1$ pages, i.\ e.\ 
$n-1$ copies of the upper half-plane, glued along $\RR$, such that the 
first edge, $\{i_1i_2\}$, of each connected component is drawn on the first page, the 
second edge, $\{i_2i_3\}$, on the second page, and so on, such that 
$i_1<i_2<\dots<i_n$, without intersections.
\end{defn}

\begin{lemma}\label{lem4}
The number of $n$-wave graphs with $m$ vertices is nonzero iff $m=kn$ for some 
integer $k$, in which case it equals \eqref{7}.
\end{lemma}
\begin{proof}
We'll construct a bijection between the set of standard rectangular tableaux of 
shape $n^k$ and the set of $n$-wave graphs with $m=kn$ vertices. More precisely, 
we'll construct a bijections between each of the mentioned sets and the set of 
the lattice words of length $m$, containing 
equally $k$ of each letter from the alphabeth $A_n=\{1, 2, \dots, n\}$.  
{\em A lattice word} means that reading it from the beginning to any letter, 
one meets at least as many 1's as 2's, at least as many 2's as 3's and so on. 

Consider a standard tableau as a notation of a Tetris game described above, 
consisting of sequential dropping of $m$ $1\times 1$ blocks. Writing sequentially 
for each dropping block the number of the column where it drops, one gets a word 
of length $m$ in the alphabeth $A_n$. It is a lattice word, because 
we require that the blocks form a Young diagram on every step. For a rectangular 
tableau of shape $n^k$ we'll get exactly $k$ copies of each letter. Conversely, each 
word of length $m$ in the alphabeth $A_n$ can be considered as a notation of some 
Tetris game consisting of sequential dropping of $m$ $1\times 1$ blocks; each 
letter shows where the corresponding block drops. Lattice words restrict us to 
the games in which dropped blocks form a Young diagram on every step, 
and words containing $k$ copies of each letter, correspond to the tableaux of 
rectangular $n^k$ shape. The constructed mappings are mutually inverse. Thus, we 
got a bijection.

Now we'll construct another one. First we'll work with case $n=2$. 
Use induction on $k$. For $k=0$ each of the sets contain just one element: empty
word and empty graph---these sets are bijective. 
Suppose that we got a bijection for all the values of $k$ less than given. 
Let $\alpha_1\dots\alpha_{2k}$ be a lattice word 
containing $k$ 1's and $k$ 2's. By definition, $\alpha_1=1$. Let $i$ be the 
smallest index such that the subword $\alpha_1 \dots \alpha_i$ contains the 
equal number of 1's and 2's. Then we draw an edge from $1$ to $i$ on the upper 
half-plane. Deleting $1$ and $i$ from the given word, we obtain either one or 
two words. They are lattice words of length less than $2k$, containing equal 
number of 1's and 2's each. By induction, 
we know the outerplanar graphs corresponding to them. Draw the graph 
corresponding to the word $\alpha_2 \dots \alpha_{i-1}$, on the vertices 
$2, \dots, i-1$, and the graph corresponding to the word $\alpha_{i+1} \dots 
\alpha_{2k}$, on the vertices ${i+1, \dots, 2k}$; instead of the standard 
indexing of vertices $\{1,2,\dots,\}$ of outerplanar graphs, preserving only 
the order; and we obtain an outerplanar graph with vertices $1, 2,\dots, 2k$. 
Conversely, setting for each edge $(i,j)$ of an outerplanar graph with $2k$ 
vertices, with $i<j$, $\alpha_i=1$ and $\alpha_j=2$, we obtain a word 
$\alpha_1 \dots \alpha_{2k}$ with all the necessary conditions. One can check, 
inductively, that the constructed mappings are mutually inverse. Thus, we got 
a bijection for the case $n=2$. 

In the general case, for $n>2$, for each lattice word with equally $k$ 1's, 2's, 
\ldots, and $n$'s; its subwords containing all of the two subsequent letters: 
1 and 2, or 2 and 3, \ldots, or $n-1$ and $n$, are lattice words in two-letter 
alphabets with equally $k$ copies of each letter. Drawing the outerplanar 
graph corresponding to the first subword on the first page, to the  
second subword---on the second page, and so on, with the numeration of vertices, 
coming from the original lattice word using $n$ letters, we obtain a wave graph. 
Conversely, setting for each edge $(i, j)$ lying on $N$-th page of a wave graph, 
where $1\leq N\leq n-1$, with $i<j$, $\alpha_i=N, \alpha_j=N+1$, we obtain a 
lattise word $\alpha_1 \dots \alpha_m$ having $k$ 1's, $k$ 2's, \ldots, $k$ 
$n$'s. The constructed mappings are mutually inverse. Thus we constructed a 
bijection. 

To complete the proof of the lemma, the same as for the proof of Lemma \ref{lem3}, 
one can use the hook length formula \cite{FRT}.
\end{proof}

\begin{defn}
For an $n$-wave graph $G$, denote ${\cal O}(G)$ the set of directed  
$n$-hypergraphs, the underlying wave graph of which is $G$. 
A directed $n$-hypergraph means 
a set of vertices $V$ equipped with a set of directed $n$-edges, i.\ e.\ 
totally ordered subsets of $n$ elements of $V$. For $V=\{1,2,\dots, m\}$, 
the underlying wave graph is a simple graph with   
the set of vertices $V$, constructed from the given directed $n$-hypergraph
by drawing for each of its edges $(i_1,\dots,i_n)$ a path $(j_1,j_2), (j_2,j_3), 
\dots, (j_n-1,j_n)$ where $(j_1,\dots,j_n)$ is such permutation of 
$(i_1,\dots,i_n)$ that $j_1<j_2<\dots<j_n$ and we draw each edge $(j_N,j_{N+1})$ 
on the $N$-th page of our book. 
\end{defn}

To define an element of ${\cal O}(G)$, we have to introduce a linear ordering on 
every connected component of $G$. Hence, the number of elements of ${\cal O}(G)$ 
is equal to $(n!)^k$ for a $n$-wave graph with $m=kn$ vertices.

\begin{defn}For $g\in{\cal O}(G)$ denote 
\begin{equation}\label{8}
b_g=x_1(g)\otimes \dots \otimes x_m(g)\in V^{\otimes m} ,
\end{equation} 
where $x_i(g)=x_{\orn i}$, and $\orn i$ is the ordinal number of vertex $i$ in the 
unique $n$-edge of $g$ it belongs. Also denote 
\begin{equation}\label{9}
t_G=\sum_{g\in{\cal O}(G)}(-1)^{\inv g} b_g ,
\end{equation}
where $\inv g$ is the number of inversions in $g$, i.\ e.\ the number of pairs of 
vertices $i<j$ ordered in an opposite way in a directed $n$-edge of $g$.  
\end{defn} 

Note that $t_G$ is equal to the tensor product of wedge products 
$x_1\wedge\dots\wedge x_n$ corresponding to the connected components, i.\ e.\ 
waves, of $G$. In more detail, for a graph $G$ with waves $(1,2,\dots,n), 
(n+1,\dots, 2n), \dots, ((k-1)n+1, \dots, kn)$
\begin{equation}\label{10}
t_G=(x_1\wedge\dots\wedge x_n)^{\otimes k},
\end{equation} 
and for other $n$-wave graphs, $t_G$ can be obtained from \eqref{10} by 
corresponding permutation (not unique, in general).

\begin{theorem}
Tensors $t_G$ parametrized by all $n$-wave graphs with $m$ vertices, form a basis 
in the space of $SL(n)$-invariants in $V^{\otimes m}$.
\end{theorem}
\begin{proof}
Lemmas \ref{lem3} and \ref{lem4} show that the number of $n$-wave graphs with $m$ 
vertices is exactly the same as the dimension of the corresponding space of 
$SL(n)$-invariants. The paragraph just before the theorem, Lemma \ref{lem2} and 
the fact that a tensor product of invariants is invariant and permutations 
of the components of a tensor product commute with the action of $SL(n)$, show 
that $t_G$ are $SL(n)$-invariant. 
Hence if we prove linearly independence of the set of $t_G$, 
our theorem will be proven. 
The proof is completely analogous to the proof of the particular case $n=2$ given 
in my articles \cite{M1,M2}. 

Denote $B$ the standard basis of $V^{\otimes m}$, consisting of $n^m$ tensor 
products $X_1\otimes\dots\otimes X_m$ with $X_1, \dots, X_m \in \{x_1, \dots, 
x_n\}$. We suppose that $B$ is ordered lexicographically. Note that $b_g\in B$ for 
any $g\in {\cal O}(G)$. Notice that for each $G$ exists exactly one $g_0\in 
{\cal O}(G)$ without inversions---with the orientation of each $n$-edge from the 
left to the right. Changing of orientation of $g_0$ increases $b_G$ in the 
lexicographical order of $B$. It means that $b_{g_0}$ is the minimal element with 
a non-zero coefficient in the decomposition of $t_G$ in the basis $B$. Denote 
$b_G=b_{g_0}$. Notice that after deleting all $x$'s and all $\otimes$ signs 
in $b_G$, the remaining indexes form the lattice word corresponding to $G$ by 
the bijection constructed in the proof of Lemma \ref{lem4}. So, we have 
$f^{n^k}$ elements $b_G$ ---one for each $G$. 

To prove the linear independence of the set of $t_G$, we can show that the rank 
of the $f^{n^k} \times n^m$ matrix of the coefficients of $t_G$ in the basis $B$ is equal 
to $f^{n^k}$. To do that, we can find a non-zero $f^{n^k} \times f^{n^k}$ minor of that matrix. 
Consider the $f^{n^k} \times f^{n^k}$ submatrix with rows numbered by $G$ ordered the same way 
as $b_G$, and columns corresponding to $b_G$. As we noticed above,  
$b_G$ is the first element with a nonzero coefficient in the row $G$, and  
this coefficient equals 1 by definition. So, this matrix is unipotent, its determinant 
equals $1$, that completes the proof of the linear independence of $t_G$.
\end{proof}

\begin{acn}
I would like to thank my thesis formal and informal advisors, Fan Chung Graham 
and Alexandre Kirillov, the supervisor of the 
part of the research supported by ONR grant, Andre Scedrov, 
our Mathematics Department and Graduate Group Chairpersons, 
Dennis DeTurck and Chris Croke, as well as UPenn Professors Christos 
Athanasiadis, Ching-Li Chai,
Ted Chinburg, Murray Gerstenhaber, Herman Gluck, Michael Larsen, 
David Shale, Stephen Shatz, Herb Wilf, and Wolfgang Ziller for
useful discussions, 
and my Gorgeous and Brilliant Wife, Bette, for her total support and love.
\end{acn}

\end{document}